\documentclass[12pt]{amsart}
\usepackage{amsmath,amscd,amssymb,amsfonts,graphics}
\newcommand{\h}{\hbox}
\newcommand{\q}{\quad}

\newcommand{\bs}{\par\bigskip}
\newcommand{\ms}{\par\medskip}
\newcommand{\sk}{\par\smallskip}
\newcommand{\bsn}{\par\bigskip\noindent}
\newcommand{\msn}{\par\medskip\noindent}
\newcommand{\skn}{\par\smallskip\noindent}
\newcommand{\ges}{\geqslant}
\newcommand{\les}{\leqslant}
\newcommand{\1}{\hskip1pt}
\newcommand{\mcap}{\hbox{$\bigcap$}}
\newcommand{\mcup}{\hbox{$\bigcup$}}
\newcommand{\msqcup}{\hbox{$\bigsqcup$}}

\newcommand{\mopl}{\hbox{$\bigoplus$}}
\newcommand{\mprod}{\hbox{$\prod$}}
\newcommand{\Fb}{{\mathcal F}^{\ssb}}
\newcommand{\Fr}{{\mathcal F}_{(r)}}
\newcommand{\A}{{\mathcal A}}

\newcommand{\F}{{\mathcal F}}

\newcommand{\Hc}{{\mathcal H}}

\newcommand{\Sc}{{\mathcal S}}

\newcommand{\PP}{{\mathbb P}}
\newcommand{\Q}{{\mathbb Q}}
\newcommand{\C}{{\mathbb C}}
\newcommand{\N}{{\mathbb N}}

\newcommand{\RR}{{\mathbf R}}

\newcommand{\Z}{{\mathbb Z}}

\newcommand{\de}{\delta}
\newcommand{\ep}{\varepsilon}

\newcommand{\Co}{{}\,\overline{\!C}{}}
\newcommand{\Do}{{}\,\overline{\!D}{}}

\newcommand{\Wo}{\overline{W}}

\newcommand{\Yo}{\overline{Y}}
\newcommand{\Ct}{{}\,\widetilde{\!C}{}}
\newcommand{\Pt}{{}\,\widetilde{\!P}{}}

\newcommand{\Xt}{{}\,\widetilde{\!X}{}}
\newcommand{\Ht}{\widetilde{H}}
\newcommand{\Ff}{F_{\!f,0}}
\newcommand{\al}{\alpha}
\newcommand{\la}{\lambda}
\newcommand{\Si}{\Sigma}

\newcommand{\dd}{\partial}

\newcommand{\mm}{^{\mathfrak m}\!\1}
\newcommand{\mmp}{^{{\mathfrak m}^+}\!\!\1}
\newcommand{\Yc}{Y^{\circ}}

\newcommand{\bl}{\bigl}
\newcommand{\br}{\bigr}
\newcommand{\pl}{\1{+}\1}
\newcommand{\mi}{\1{-}\1}
\newcommand{\eq}{\,{=}\,}
\newcommand{\col}{\1{:}\1}

\newcommand{\ssb}{\raise.15ex\h{${\scriptscriptstyle\bullet}$}}
\newcommand{\ssc}{\,\raise.15ex\h{${\scriptstyle\circ}$}\,}
\newcommand{\onto}{\mathop{\rlap{$\to$}\hskip2pt\hbox{$\to$}}}
\newcommand{\into}{\hookrightarrow}
\newcommand{\simto}{\,\,\rlap{\hskip1.5mm\raise1.4mm\hbox{$\sim$}}\hbox{$\longrightarrow$}\,\,}
\newcommand{\plim}{\rlap{\raise-5.5pt\h{$\,\leftarrow$}}{\rm lim}}
\begin{document}
\title[Lowest non-zero vanishing cohomology]
{Lowest non-zero vanishing cohomology\\of holomorphic functions}
\author[M. Saito]{Morihiko Saito}
\address{M. Saito : RIMS Kyoto University, Kyoto 606-8502 Japan}
\email{msaito@kurims.kyoto-u.ac.jp}
\begin{abstract} We study the vanishing cycle complex $\varphi_fA_X$ for a holomorphic function $f$ on a reduced complex analytic space $X$ with $A$ a Dedekind domain (for instance, a localization of the ring of integers of a cyclotomic field, where the monodromy eigenvalue decomposition may hold after a localization of $A$). Assuming the perversity of the shifted constant sheaf $A_X[d_X]$, we show that the lowest possibly-non-zero vanishing cohomology at $0\in X$ can be calculated by the restriction of $\varphi_fA_X$ to an appropriate nearby curve in the singular locus $Y$ of $f$, which is given by intersecting $Y$ with the intersection of sufficiently general hyperplanes in the ambient space passing sufficiently near 0. The proof uses a Lefschetz type theorem for local fundamental groups. In the homogeneous polynomial case, a similar assertion follows from Artin's vanishing theorem. By a related argument we can show the vanishing of the non-unipotent monodromy part of the first Milnor cohomology for many central hyperplane arrangements with ambient dimension at least 4.
\end{abstract}
\maketitle
\centerline{\bf Introduction}
\bsn
Let $f$ be a holomorphic function on a reduced complex analytic space $X$ of pure dimension $d_X=n{+}1$.
We are interested in the nearby and vanishing cycle complexes $\psi_fA_X[n]$, $\varphi_fA_X[n]$ with $A$ a field or more generally a Dedekind domain, for instance, $A=\Z$, $\Q$, $\C$. If $A$ contains a primitive $m$-th root of unity with $m$ the smallest positive integer such that $T_s^m=1$ (where $T_s$ is the semi-simple part of the Jordan decomposition of the monodromy $T$ defined after a localization of $A$), then we have the monodromy eigenvalue decompositions of $\psi_fA_X[n]$, $\varphi_fA_X[n]$ (after a localization of $A$), see (1.3) below.
\sk
Assume the perversity of the shifted constant sheaf $A_X[d_X]$ is satisfied (for instance, $X$ is smooth), see \cite{BBD}, etc. This is a fairly strong condition; $X$ must be pure-dimensional, and the intersection of irreducible components must have codimension 1 if $X$ has two irreducible components. (This follows from Proposition~1 below.) Note also that the affine cone of a smooth projective variety $Y$ cannot satisfy the perversity condition if $Y$ has non-trivial primitive cohomology of degree $j\in[1,d_Y{-}1]$, for instance, an abelian variety of dimension at least 2
(since the stalk at the origin of its intersection complex is given by the primitive cohomology of $Y$ using the Thom-Gysin sequence).
\sk
The perversity of $A_X[d_X]$ implies that of $\psi_fA_X[n]$, $\varphi_fA_X[n]$ with $n=d_X\mi 1$, since this property is preserved by the shifted nearby and vanishing cycle functors $^p\psi_f:=\psi_f[-1]$, $^p\varphi_f:=\varphi_f[-1]$, see \cite{KS}, \cite{Sc} (and \cite{Br} for the case $A$ is a field). This imposes a quite strong condition on $\varphi_fA_X[n]$ as is explained in \cite{DS1}. We have for instance
$$\F_{(j)}:=\Hc^{-j}(\varphi_fA_X[n])=0\q\h{if}\q j>r:=\dim{\rm Supp}\,\varphi_fA_X.
\leqno(1)$$
This follows from the stability by shifted nearby and vanishing cycle functors explained above. Indeed, the pull-back functor for $\{0\}\into X$ is expressed by the iteration of the mapping cones of nearby and vanishing cycle functors associated with holomorphic functions such that the intersection of their zero loci is 0. Here we can apply also the following.
\bsn
{\bf Proposition~1.} {\it Assume $\Fb\in D^b_c(Y,A)$ satisfies the semi-perversity $\mm D^{\ges 0}$ $($see $(1.1)$ below$)$, where $Y={\rm Supp}\,\Fb$ is a reduced complex analytic space of dimension $r$. Let $\Sc$ be a Whitney stratification of $Y$ compatible with $\Fb$. Let $Y^i$ be the union of strata $S\in\Sc$ with $\dim S\ges i$. Let $j_i:Y^i\into Y$ be the canonical inclusions.
Then, for $i\in[1,r]$, there are canonical isomorphisms
$$\Hc^{-k}\1\Fb\simto\Hc^{-k}\1\RR(j_i)_*(\Fb|_{Y^i})\q\q(\1k>i\1),
\leqno(2)$$
and the inclusion}
$$\Hc^{-i}\1\Fb\into\Hc^{-i}\1\RR(j_i)_*(\Fb|_{Y^i}).
\leqno(3)$$
\ms
This is quite useful, although its proof is an easy exercise of derived categories, see (1.2) below. In the case $\Fb=\varphi_fA_X[n]$ with $A=\Q$ or $\C$, the assertions~(2--3) for $i\eq 1$, and the assertions (2) for $(k,i)=(r,r{-}1)$ and (3) for $i\eq r$ were proved and announced in \cite[Theorem~0.1 and (3.5.1--2)]{DS1} respectively. Here we use also the Leray type spectral sequence and (1). In (3.5.1--2) we can consider only the monodromy eigenvalues $\la$ with $\dim{\rm Supp}\,\varphi_{f,\la}\C_X=r$, and take the direct sum of $\varphi_{f,\la}\C_X$ over such $\la$, since the other $\la$ can be neglected by (1). It is not difficult to replace $\Q$ or $\C$ with $A$ using \cite{KS}, \cite{Sc}.
See also \cite[Theorem~3.1]{MPT} where it is assumed that $\Fb=\varphi_f\Z_X[n]$ and moreover $i=k{-}1$ in (2).
It does not seem completely trivial to deduce the general case from it, unless the {\it sheaf version\1} as in Proposition~1 is employed, since {\it hypercohomology\1} is used there so that only the assertion for each stalk is proved essentially. (Note that the {\it topological cone theorem\1} is needed to show the assertion on hypercohomology.)
It may be viewed as a minimal extension of \cite[(3.5.1--2)]{DS1}, since we get the same $i$ restricting to the lowest degree part $k\eq r$.
\sk
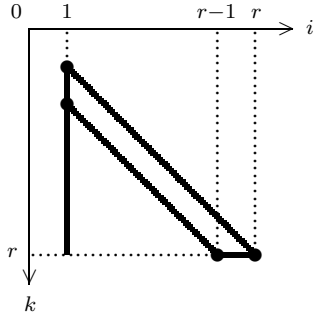
$$\q\h{$\setlength{\unitlength}{5mm}
\begin{picture}(8,8)
\linethickness{.1mm}
\multiput(1,7.4)(0,-0.2){5}{\circle*{.1}}
\multiput(0.1,1.5)(0.2,0){30}{\circle*{.1}}
\multiput(6,7.4)(0,-0.2){30}{\circle*{.1}}
\multiput(5,7.4)(0,-0.2){30}{\circle*{.1}}
\put(0,7.5){\line(1,0){7}}
\put(0,.7){\line(0,1){6.8}}
\put(-.5,7.8){$\scriptstyle 0$}
\put(-.6,1.45){$\scriptstyle r$}
\put(5.9,7.8){$\scriptstyle r$}
\put(4.45,7.8){$\scriptstyle r{-}1$}
\put(.85,7.8){$\scriptstyle 1$}
\put(7,7.5){\raise-.635mm\h{\!\!\!$\scriptstyle >\,\,\,i$}}
\put(6,1.5){\circle*{.35}}
\put(5,1.5){\circle*{.35}}
\put(1,5.5){\circle*{.35}}
\put(1,6.5){\circle*{.35}}
\put(-.18,.7){\rotatebox{90}{$\scriptstyle<$}}
\put(-.14,0){$\scriptstyle k$}
\linethickness{.6mm}
\put(1,1.5){\line(0,1){5}}
\put(5,1.5){\line(1,0){1}}
\qbezier(1,5.5)(3,3.5)(5,1.5)
\qbezier(1,6.5)(3.5,4)(6,1.5)
\end{picture}$}$$
\sk
The proof of Proposition~1 is essentially the same as the argument suggested in {\it loc.\,cit.,} where it is stated as follows: ``Similarly (3.5.2) follows from Theorem~0.1 by induction on stratum". If an assertion concerning a Zariski-open immersion is proved ``by induction on stratum", it usually means decomposing the open immersion so that the assertion is reduced to the case where the complement of an open immersion is a stratum, or more generally, a union of strata of the same dimension.
\sk
We apply Proposition~1 to the case
$$Y:={\rm Supp}\,\Fb\q\h{with}\q\Fb:=\varphi_fA_X[n].$$
In this paper, the singular locus of $f$ is defined to be the support of $\varphi_fA_X[n]$.
(It does not seem very clear whether it can be defined by using a stratification of $X$.)
Let $\Ff$ denote the Milnor fiber of $f$ at $0\in Y$. By \cite{Hi} we have the isomorphisms
$$\Ht^{n-j}(\Ff,A)=\F_{(j),0}\q\h{with}\q\F_{(j)}:=\Hc^{-j}\Fb.
\leqno(4)$$
This holds even if $\Ff=\emptyset$, where the reduced cohomology is defined to be the cohomology of the mapping cone of $A\to 0$.
Here note that
$$\F_{(j)}=\Hc^{-j}(\psi_fA_X[n])\q\h{if}\q j\les r\les n{-}1.$$
\sk
Assume $(Y,0)$ is a closed analytic subset of $(\C^N,0)$. Let $B_{\ep}$ be an open ball of radius $\ep\ll 1$ in $\C^N$. Let $p:(Y,0)\to(\C^{r-1},0)$ be the restriction of a sufficiently general linear projection $p_{\C^N}:(\C^N,0)\to(\C^{r-1},0)$.
\msn
{\bf Theorem~1.} {\it Assume the perversity of $A_X[d_X]$ is satisfied. Then we have the isomorphism
$$\Ht^{n-r}(\Ff,A)=\Gamma(C_{\ep},\Fr|_{C_{\ep}}),
\leqno(5)$$
with $C_{\ep}\subset Y$ a curve defined by $B_{\ep}\cap p^{-1}(z)$. Here $z$ is a sufficiently general point of an open ball $B'_{\de}\subset\C^{r-1}$ with radius $\de\ll\ep\ll 1$.}
\ms
If $r\les n{-}1$, we have $\Fr|_{C_{\ep}}=\Hc^{n-r}\psi_{f_W}A_W$ with $W:=B_{\ep}\cap p_X^{-1}(z)$, $f_W:=f|_W$, and $p_X$ the restriction of the projection $p_{\C^N}$ to $X$, since $p_{\C^N}$ and $z\in B'_{\de}$ are sufficiently general.
\sk
Theorem~1 is equivalent to the following.
\msn
{\bf Theorem~2.} {\it In the assumption of Theorem~$1$, the Milnor cohomology $\Ht^{n-r}(\Ff,A)$ is isomorphic to the kernel of the canonical morphism
$$\iota_0:\Gamma\bl(C^{\circ}_{\ep},\Fr|_{C^{\circ}_{\ep}}\br)\to\Gamma\bl(D,\bl((j_C)_*(\Fr|_{C^{\circ}_{\ep}})/(\Fr|_{C_{\ep}})\br)|_D\br).
\leqno(6)$$
Here $C^{\circ}_{\ep}:=C_{\ep}\cap Y^r$, $D:=C_{\ep}\setminus C^{\circ}_{\ep}$ with inclusions $j_C:C^{\circ}_{\ep}\into C_{\ep}$, $i_D:D\into C_{\ep}$.}
\ms
Theorems~1 and 2 are equivalent to each other, since $\Gamma(C_{\ep},\Fr|_{C_{\ep}})$ is given by ${\rm Ker}\,\iota_0$ (using the $A$-coefficient version of \cite[Theorem~0.1]{DS1}).
The target of $\iota_0$ is viewed as the {\it obstruction\1} for an section of $\Fr|_{C^{\circ}_{\ep}}$ to be extended to that of $\Fr|_{C_{\ep}}$.
The sheaf whose global sections form the target of $\iota_0$ is called the {\it obstruction local system}.
This is a generalization of the short exact sequence in \cite[Theorem~0.1]{DS1} (where $K_x$ is the obstruction), and is closely related to an inequality in \cite[Theorem~3.4\1(b)]{MPT}.
\sk
It is also possible to consider the canonical morphisms
$$\aligned\iota_1:\Gamma\bl(C^{\circ}_{\ep},\Fr|_{C^{\circ}_{\ep}}\br)\to\Gamma\bl(D,(j_C)_*(\Fr|_{C^{\circ}_{\ep}})|_D\br),\\
\iota_2:\Gamma\bl(D,\Fr|_D\br)\into\Gamma\bl(D,(j_C)_*(\Fr|_{C^{\circ}_{\ep}})|_D\br).\endaligned$$
The kernel of $\iota_0$ is isomorphic to the kernel of
$$\iota_1\mi\iota_2:\Gamma\bl(C^{\circ}_{\ep},\Fr|_{C^{\circ}_{\ep}}\br)\oplus\Gamma\bl(D,\Fr|_D\br)\to\Gamma\bl(D,(j_C)_*(\Fr|_{C^{\circ}_{\ep}})|_D\br),
\leqno(7)$$
by gluing sections. This is useful {\it only\1} to see the relation with a formulation in \cite[Theorem~5.5]{MPT} (which was seriously misstated in the first two versions).
Note that $\iota_2$ is {\it injective\1} by the $A$-coefficient version of \cite[Theorem~0.1]{DS1}.
This implies that ${\rm Ker}(\iota_1\mi\iota_2)$ is uniquely determined by its projection to the first direct factor of the source of (7).
\sk
By the $A$-coefficient version of \cite[(3.5.2)]{DS1} with monodromy eigenvalues forgotten, we have the isomorphism
$$\Ht^{n-r}(\Ff,A)=\Gamma\bl(Y^{r-1}_{\ep,\de},\Fr|_{Y^{r-1}_{\ep,\de}}\br)\q(0<\de\ll\ep\ll 1),
\leqno(8)$$
where
$$Y^k_{\ep,\de}:=Y^k\cap B_{\ep}\cap p^{-1}(B''_{\de})\q\q(k\in\N),$$
with $Y^k$ as in Proposition~1 and $B''_{\de}\subset\C^{r-1}$ an open ball of radius $\de$.
Indeed, the right-hand side is independent of $0<\de\ll\ep\ll 1$, and we have a constant inductive system. Set
$$Z_{\ep,\de}:=Y^{r-1}_{\ep,\de}\setminus Y^r_{\ep,\de},$$
with $i_Z:Z_{\ep,\de}\into Y^{r-1}_{\ep,\de}$, $j_Y:Y^r_{\ep,\de}\into Y^{r-1}_{\ep,\de}$ natural inclusions.
\sk
The right-hand side of (8) is identified with the kernel of
$$\iota'_0:\Gamma\bl(Y^r_{\ep,\de},\Fr|_{Y^r_{\ep,\de}}\br)\to\Gamma\bl(Z_{\ep,\de},\bl((j_Y)_*(\Fr|_{Y^r_{\ep,\de}})/(\Fr|_{Y^{r-1}_{\ep,\de}})\br)|_{Z_{\ep,\de}}\br).
\leqno(9)$$
Here the injectivity of $\Fr|_{Y^{r-1}_{\ep,\de}}\to(j_Y)_*(\Fr|_{Y^r_{\ep,\de}})$ follows from the $A$-coefficient version of \cite[Theorem 0.1]{DS1}.
This identification means that the target of $\iota'_0$ is the {\it obstruction\1} for a section of $\Fr|_{Y^r_{\ep,\de}}$ to be extended over $Y^{r-1}_{\ep,\de}$.
(This is closely related to the inequality in \cite[Theorem 3.4\,(b)]{MPT} via an exact sequence in \cite[Theorem 0.1]{DS1}.)
\sk
We have a Lefschetz type theorem for local fundamental groups, see Theorem~(2.1) below. This is proved by applying the long exact sequence of homotopy groups associated with a topologically fibered space. Using this, we can deduce that the above kernel does not change by restricting to the curve $C_{\ep}$. So Theorems~1 and 2 follow.
\sk
An assertion related to Theorem~2 seems to have been studied in \cite[Theorem 5.5]{MPT} using topological CW-complex models due to Siersma and Tib\u ar when $r\eq 2$. The assertion, however, seems to have been misstated rather seriously in the first version, and corrected in the third version following some suggestion. (This kind of misstatement is not called ``typo" usually.)
Here it seems more natural to use the {\it difference kernel,} rather than the {\it intersection of images,} in order to glue sections as is explained after (7). This intersection may come from the CW complex model argument. Before Lemma~5.2, {\it loc.\,cit.,} it does not seem very clear whether it is conjectured that the source of the morphism (6) is isomorphic to that of (9) in the case $r=2$. Indeed, it is written before it as follows: ``Let us point out that in general this representation cannot be related by a Lefschetz slicing argument to the monodromy representation defined at (3.19)." At least it does not seem to be mentioned that this isomorphism could imply immediately a variant of Theorem~5.5, {\it loc.\,cit}. It is moreover stated at the end of the introduction (even in the third version) as follows: ``Based on it, we obtain in Corollary~5.3 and Theorem~5.5 a description of $H^{n-s}(F)$ in terms of the invariant submodule of a monodromy representation (5.6) on the transversal Milnor fibre of the $s$-dimensional strata, which is in fact a genuine vertical monodromy representation. The similarity with the terms and results of Theorem~3.4 is striking. But surprisingly, this monodromy representation is totally different from (1.1) used in Theorem~3.4, and in general cannot be deduced from it by a Lefschetz slicing argument." If this is true, the suggestion mentioned above may lose its ground.
\sk
Our formulation makes the reduction to the case $r\eq 2$ unnecessary, and may be viewed actually as a reduction to the case $r\eq 1$. In the reduced homogeneous polynomial case, we have the following.
\msn
{\bf Corollary~1.} {\it Assume $X=\C^{n+1}$, and $f$ is a homogeneous polynomial. Let $W\subset X$ be a sufficiently general $(n{-}r{+}2)$-dimensional affine subspace. Then we have a canonical isomorphism compatible with mixed Hodge structure}
$$H^{n-r}(\Ff)=\Gamma(W,\Hc^{n-r}(\varphi_{f_W}A_W)).
\leqno(10)$$
\ms
Note that $r\,{<}\,n$ if $f$ is reduced. This corollary follows also from the {\it weak Lefschetz type theorem\1} (or Artin's vanishing theorem \cite{BBD}), see Remark~(2.3) below. Corollary~1 says that we do not have to take the intersection with an open ball $B_{\ep}$ in the case of reduced homogeneous polynomials. For central reduced hyperplane arrangements, it means that the first Milnor cohomology is calculated by the global sections of the first vanishing cohomology sheaf $\Hc^1(\psi_{f_W}A_W)$ for a {\it non-central\1} reduced plane arrangement defined by $f_W$ in $W\cong\C^3$. As a corollary, we get a simple proof of the vanishing of $H^1(\Ff)_{\la}$ with $\la=\exp(\pm 2\pi i/6)$ in the case of a reflection hyperplane arrangement of type $G_{31}$, see \cite{BDY}, \cite{ac}. It turns out, however, that this is actually a corollary of \cite[(3.5.1--2)]{DS1}, see Corollary~(2.4) below. Theorem~(2.4) below implies that hyperplane arrangements in $\C^{n+1}$ ($n\ges 3$) cannot have non-vanishing $H^1(\Ff)_{\la}$ ($\la\,{\ne}\,1$) except for some quite special cases, see (2.4), (3.5) below.
\sk
This work was partially supported by JSPS Kakenhi 15K04816.
I thank A. Dimca and J.~Sch\"urmann for useful information about \cite{Di1}, \cite{Sc}.
\sk
In Section~1 we recall the definition of semi-perversities, and prove Proposition~1.
In Section~2 we prove a Lefschetz type theorem for local fundamental groups, which implies Theorems~1 and 2.
In Section~3 we give some example.
\bs\bs
\vbox{\centerline{\bf 1. Preliminaries}
\bsn
In this section we recall the definition of semi-perversities, and prove Proposition~1.}
\msn
{\bf 1.1.~Semi-perversity conditions.} Let $X$ be a reduced complex analytic space. Let $A$ be a noetherian commutative ring of finite global dimension. We have the following semi-perversity conditions for $\Fb\in D^b_c(X,A)$\,:
$$\aligned\mm D^{\les j}&:\,\Hc^ki_S^*\Fb=0\q\q(k>j\mi d_S),\\
\mm D^{\ges j}&:\,\Hc^ki_S^!\Fb=0\q\q(k<j\mi d_S),\endaligned
\leqno(1.1.1)$$
where $S$ runs over any strata of a Whitney stratification with $i_S$ the inclusion of $S$ and $d_S=\dim S$, see \cite{BBD}, etc.
(Here $\mm$ stands for the {\it middle\1} perversity.)
We say that $\Fb$ satisfies the perversity condition if these two semi-perversity conditions hold for $j=0$.
Note that these are not dual to each other unless $A$ is a field.
The stability of these semi-perversity conditions by shifted nearby and vanishing cycle functors is shown in \cite{KS}, \cite{Sc}.
\sk
In the case $A$ is a Dedekind domain, we have the following dual semi-perversity conditions:
$$\aligned\mmp D^{\les j}&:\,\Hc^ki_S^*\Fb=\begin{cases}0&(k>j\mi d_S\pl 1),\\\h{torsion}&(k=j\mi d_S\pl 1),\end{cases}\\
\mmp D^{\ges j}&:\,\Hc^ki_S^!\Fb=\begin{cases}0&(k<j\mi d_S),\\\h{torsion-free}&(k=j\mi d_S).\end{cases}\endaligned
\leqno(1.1.2)$$
These are the dual of $\mm D^{\ges -j}$ and $\mm D^{\les -j}$ respectively, see \cite[3.3]{BBD}.
We say that the {\it strong perversity\1} holds if $\mm D^{\les 0}$ and $\mmp D^{\ges 0}$ are satisfied.
The stability of these semi-perversities by shifted nearby and vanishing cycle functors seems to be shown in \cite{Sc}. (Here the structure of arguments seems rather complicated; a very general assertion is stated in Theorem 6.0.2, and one can apply this to various situations as is explained in Example 6.0.2, where the explanation about $\mmp D^{\les j}$ seems to be skipped.)
\sk
When $A$ is the ring of integers of a number field, this stability may be related to the one in the finite field coefficient case using the short exact sequence
$$0\to I\to A\to A/I\to 0,
\leqno(1.1.3)$$
where $I\subset A$ is a maximal ideal (which is a projective $A$-module) so that $A/I$ is a finite field. It is also closely related to the duality of nearby and vanishing cycle functors, see Appendix below.
\msn
{\bf Remark~1.1}\,(i). A Noetherian commutative ring $A$ is called a Dedekind domain if it is a 1-dimensional normal domain. Here every non-zero prime ideal is a maximal ideal, and the localization of $A$ by the complement of a maximal ideal (that is, the stalk of the structure sheaf of ${\rm Spec}\,A$) is a discrete valuation ring. So torsion-freeness is equivalent to projectivity for finite $A$-modules, and higher extension groups ${\rm Ext}^i$ vanish for $i>1$. Moreover any finite $A$-module is a direct sum of a torsion-free $A$-module and a torsion $A$-module.
\sk
It is well known that a localization of the ring of integers of a number field is a Dedekind domain. The reader may assume that a Dedekind domain is such a ring in this paper.
\msn
{\bf Remark~1.1}\,(ii). Assume the perversity and the strong perversity are satisfied for $A_X[d_X]$ and $A_{X\setminus Z}[d_X]$ respectively with $Z\subset{\rm Sing}\,f$ a closed analytic subset satisfying the condition $\dim Z\,{<}\,r\eq\dim{\rm Sing}\,f$, where $A$ is a Dedekind domain. Then $\Hc^{n-r}\varphi_fA_X$ is torsion-free.
\msn
{\bf 1.2.~Proof of Proposition~1.}
Let $j'_i:Y^{i+1}\into Y^i$ be natural inclusions. Set $\Si_i:=Y^i\setminus Y^{i+1}$.
There are distinguished triangles
$$\RR\Gamma_{\Si_i}(\Fb|_{Y^i})\to\Fb|_{Y^i}\to\RR(j'_i)_*(\Fb|_{Y^{i+1}})\buildrel{+1}\over\to.
\leqno(1.2.1)$$
The assumption on the semi-perversity of $\Fb$ implies that
$$\Hc^{-k}\1\RR\Gamma_{\Si_i}(\Fb|_{Y^i})=0\q\q(\1k>i=\dim\Si_i\1),
\leqno(1.2.2)$$
see (1.1). Note that $\RR\Gamma_{\Si_i}=(i_{\Si_i})_*i_{\Si_i}^!$ with $i_{\Si_i}:\Si_i\into Y^i$ the inclusion.
\sk
Since $j_i\ssc j'_i=j_{i+1}$, we get the distinguished triangles
$$\RR(j_i)_*\RR\Gamma_{\Si_i}(\Fb|_{Y^i})\to\RR(j_i)_*(\Fb|_{Y^i})\to\RR(j_{i+1})_*(\Fb|_{Y^{i+1}})\buildrel{+1}\over\to,
\leqno(1.2.3)$$
and also the vanishing:
$$\Hc^{-k}\1\RR(j_i)_*\RR\Gamma_{\Si_i}(\Fb|_{Y^i})=0\q\q(\1k>i\1).
\leqno(1.2.4)$$
Indeed, $\RR(j_i)_*$ is the right derived functor of the left exact functor $(j_i)_*$.
\sk
These imply the isomorphisms
$$\Hc^{-k}\1\RR(j_i)_*(\Fb|_{Y^i})\simto \Hc^{-k}\1\RR(j_{i+1})_*(\Fb|_{Y^{i+1}})\q\q(\1k>i\pl 1\1),
\leqno(1.2.5)$$
and the inclusions
$$\Hc^{-i-1}\1\RR(j_i)_*(\Fb|_{Y^i})\into \Hc^{-i-1}\1\RR(j_{i+1})_*(\Fb|_{Y^{i+1}}).
\leqno(1.2.6)$$
\sk
So we get the isomorphisms (2) and the inclusions (3). Hence Proposition~1 follows.
(This argument is essentially the same as the one used in \cite[Section 3.5]{DS1}, where it is noted as ``by induction on stratum", see remarks after Proposition~1.)
\msn
{\bf 1.3. Monodromy eigenvalue decomposition.} Let $f$ be a holomorphic function on a reduced complex analytic variety $X$ such that the perversity of $A_X[d_X]$ holds. Let $T=T_sT_u$ be the Jordan decomposition with $T_s,T_u$ the semi-simple and unipotent part. This can be defined in $K$, the field of fractions of $A$.
\sk
Assume $A$ contains a primitive $m$\1th root of unity with $m$ the smallest positive integer with $T_s^m={\rm id}$ (for instance $K$ is a cyclotomic field).
There is $q\in A\setminus\{0\}$ such that $T_s$ is a polynomial of $T$ with coefficients in $A_q$ (the localization of $A$ by $q$) and moreover we have the monodromy eigenvalue decompositions\,:
$$\aligned\psi_fA_{q,X}[n]=\mopl_{\la\in\mu_m}\,\psi_{f,\la}A_{q,X}[n],\\
\varphi_fA_{q,X}[n]=\mopl_{\la\in\mu_m}\,\varphi_{f,\la}A_{q,X}[n],\endaligned
\leqno(1.3.1)$$
with $\mu_m:=\{a\in A\mid a^m=1\}$, and
$$\aligned\psi_{f,\la}A_{q,X}[n]={\rm Ker}(T_s\mi\la)\subset\psi_fA_{q,X}[n],\\ \varphi_{f,\la}A_{q,X}[n]={\rm Ker}(T_s\mi\la)\subset\varphi_fA_{q,X}[n].\endaligned
\leqno(1.3.2)$$
\msn
{\bf Remark~1.3.} We can construct the projectors to ${\rm Ker}(T_s\mi\la)$ using a minimal polynomial $P(x)\in A_q[x]$ of $T$. Indeed, if there is a factorization $P(x)=P_1(x)P_2(x)$ such that $P_1(x)$, $P_2(x)$ are mutually prime in $K[x]$ (where $K$ is the filed of fractions of $A$), then there are $Q_1(x)$, $Q_2(x)\in A_q[t]$ such that
$$P_1(x)Q_1(x)+P_2(x)Q_2(x)=1\q\h{in}\,\,\,A_q[x],
\leqno(1.3.3)$$
replacing $q$ if necessary. Since $P_1(T)P_2(T)\eq 0$, we get projectors $\pi_i:=P_i(T)Q_i(T)$ ($i\eq 1,2$) such that $\pi_i^2\eq\pi_i$ and $\pi_i\ssc\pi_j\eq 0$ ($i\,{\ne}\,j$).
We can continue this if there is a factorization of $P_1(x)$ or $P_2(x)$ as above.
\bs\bs
\vbox{\centerline{\bf 2. Lefschetz type theorem}
\bsn
In this section we prove a Lefschetz type theorem for local fundamental groups, which implies Theorems~1 and 2.}
\msn
{\bf 2.1. Lefschetz type theorem for local fundamental groups.} Let $(Y,0)$ be a germ of a reduced complex analytic space of pure dimension $r$. Let $(Z,0)\subset(Y,0)$ be a closed analytic subset of dimension strictly smaller than $r$ such that $\Yc:=Y\setminus Z$ is smooth.
Let $\{Y_{(i)}\}_{i\in\N}$ be a fundamental system of open neighborhoods of $0\in Y$ satisfying the following condition:
$$\h{We have the bijection $\pi_1\bl(\Yc_{(i+1)}\br)\simto\pi_1\bl(\Yc_{(i)}\br)$ for any $i\ges i_0$,}
\leqno(2.1.1)$$
for some $i_0\in\N$, where $\Yc_{(i)}:=\Yc\cap Y_{(i)}$.
For instance, $Y_{(i)}=Y\cap B_{\ep_i}$ with $B_{\ep_i}$ an open ball with radius $\ep_i$ in an ambient space. Here the $\ep_i$ are strictly decreasing, and $\lim_{i\to\infty}\ep_i=0$.
(One may have to use closed balls (or collar neighborhoods) if one wants to apply the notion of ``deformation retract".)
\sk
We define the {\it local fundamental group\1} of $\Yc$ at 0 to be the projective limit of the constant projective system $\pi_1\bl(\Yc_{(i)}\br)$ ($i\ges i_0$),
This is independent of the choice of a fundamental system of open neighborhoods satisfying (2.1.1) (using the constancy of the projective systems, since it is a ``fundamental" system of neighborhoods.)
\sk
We have a Lefschetz type theorem for local fundamental groups as follows.
\msn
{\bf Theorem~2.1.} {\it Let $p:Y\to\C^{r-1}$ be a sufficiently general projection which is the restriction of a projection of an ambient space $\C^N$. Let $B'_{\de}\subset\C^{r-1}$ be an open ball of radius $\de$. There is a closed analytic subset $(\Si,0)\subset(\C^{r-1},0)$ together with an integer $i_1\ges i_0$ such that the following holds\1$:$ for any $i\ges i_1$, there is $\de_i>0$ such that for any $z\in B'_{\de_i}\setminus \Si$, we have the surjectivity of the canonical morphism}
$$\pi_1\bl(\Yc_{(i)}\cap p^{-1}(z)\br)\onto\pi_1\bl(\Yc_{(i)}\br).
\leqno(2.1.2)$$
\msn
{\it Proof.} We first show that $Z$ can be enlarged as long as its dimension remains strictly smaller than that of $Y$. Indeed, for such $(Z',0)\supset(Z,0)$, put $Y^{\prime\circ}_{(i)}:=Y_{(i)}\setminus Z'$. We have the commutative diagram
$$\begin{array}{ccc}\pi_1\bl(Y^{\prime\circ}_{(i)}\cap p^{-1}(z)\br)&\to&\pi_1\bl(Y^{\prime\circ}_{(i)}\br)\\ \rlap{$\downarrow$}\raise.7mm\h{$\downarrow$}&\raise5mm\h{}\raise-3mm\h{}&\rlap{$\downarrow$}\raise.7mm\h{$\downarrow$}\\
\pi_1\bl(\Yc_{(i)}\cap p^{-1}(z)\br)&\to&\pi_1\bl(\Yc_{(i)}\br)\end{array}
\leqno(2.1.3)$$
where the vertical morphisms are surjective.
So the assertion is reduced to the one for $Z'$, and we may replace $Z$ with $Z'$.
\sk
We assume that $p$ is the composition of a sufficiently general projection $p':(Y,0)\to(\C^r,0)$ with $q:(\C^r,0)\to(\C^{r-1},0)$. Note that $p'$ is a finite morphism (in particular, proper). This follows from the Weierstrass preparation theorem.
\sk
Let $(\Xi,0)\subset(\C^r,0)$ be a closed analytic subset of codimension 1 such that $p'$ is locally biholomorphic over the complement of $(\Xi,0)$ and $p'(Z)\subset\Xi$.
Since $p$ is sufficiently general, we may assume that $q^{-1}(0)\cap\Xi=\{0\}$, and the restriction of $q$ to $(\Xi,0)$ is a finite morphism. Let $(\Si,0)\subset(\C^{r-1},0)$ be a closed analytic subset of codimension 1 such that the restriction of $q$ to $\Xi$ is locally biholomorphic over the complement of $(\Si,0)$.
\sk
We now replace $Z$ so that $Z=p'{}^{-1}(\Xi)$.
Using the independency of fundamental systems of open neighborhoods satisfying condition (2.1.1), we may assume that the fundamental system of open neighborhoods is given by the pull-back by $p'$ of open balls $B''_{\ep_i}\subset\C^r$. (Indeed, (2.1.1) is satisfied, since $p'$ is a finite morphism, and is locally biholomorphic over the complement of $\Xi$.)
Here we may further shrink it by taking the intersection with $p^{-1}(B'_{\de_i})$ with $\de_i>0$ very small (depending on $\ep_i$).
We may assume that the $\de_i$ are strictly decreasing.
To see that condition (2.1.1) is still satisfied after this replacement, we have to show that
$$|q|:(B''_{\ep_i}\setminus\Xi)\cap q^{-1}(B'_{\de_i}\setminus\{0\})\to (0,\de_i)\q(i\ges i_0)
\leqno(2.1.4)$$
is a topological fibration if $\de_i>0$ is very small (since $p'$ is locally biholomorphic over the complement of $\Xi$).
This can be proved by using the Thom-Mather theory, since $(\Xi,0)$ is finite over $(\C^{r-1},0)$ and $\dd B''_{\ep_1}\cap q^{-1}(z)$ is a circle for $z\in B'_{\de_i}$.
\sk
We have also a topological fibration
$$p:\Yc_{(i)}\cap p^{-1}(B'_{\de_i}\setminus\Si)\to B'_{\de_i}\setminus\Si,
\leqno(2.1.5)$$
replacing $\de_i>0$ if necessary (depending on $\ep_i$). Indeed, since $p'$ is a finite unramified covering over the complement of $(\Xi,0)$ in $(\C^r,0)$, the assertion is reduced to showing the topological fibration
$$(B''_{\ep_i}\setminus\Xi)\cap q^{-1}(B'_{\de_i}\setminus\Si)\to B'_{\de_i}\setminus\Si,
\leqno(2.1.6)$$
induced by $q$. But this follows from the relation between $\Xi$ and $\Si$ explained above. So we get the desired fibration.
\sk
By the long exact sequence of homotopy groups associated with the topological fibration, we have the exact sequence
$$\pi_1\bl(\Yc_{(i)}\cap p^{-1}(z)\br)\to\pi_1\bl(\Yc_{(i)}\cap p^{-1}(B'_{\de_i}\setminus\Si)\br)\to\pi_1(B'_{\de_i}\setminus\Si)\to 1,
\leqno(2.1.7)$$
with $z\in B'_{\de_i}\setminus\Si$. Here $\Yc_{(i)}\cap p^{-1}(z)$ is connected, since $q$ is sufficiently general (and $p=q\ssc p'$).
\sk
This exact sequence implies the surjection
$$\pi_1\bl(\Yc_{(i)}\cap p^{-1}(z)\br)\to\pi_1\bl(\Yc_{(i)}\cap p^{-1}(B'_{\de_i})\br).
\leqno(2.1.8)$$
Indeed, a loop around an irreducible component of $\Si\cap B'_{\de_i}$ can be lifted to a loop around an irreducible component of $\Yc_{(i)}\cap p^{-1}(\Si\cap B'_{\de_i})$, and its image in $\pi_1\bl(\Yc_{(i)}\cap p^{-1}(B'_{\de_i})\br)$ vanishes.
\sk
As is explained before (2.1.4), the fundamental system of open neighborhoods can be replaced by taking the intersection with $p^{-1}(B'_{\de_i})$. So Theorem~(2.1) follows.
\msn
{\bf 2.2. Proofs of Theorems~1 and 2.} It is sufficient to prove Theorem~2. By the argument after it in the introduction, it is enough to show that the kernel of $\iota'_0$ does not change by restricting to the curve $C_{\ep}=p^{-1}(z)$.
By Theorem~(2.1) the source of the morphism does not change.
(This seems to be conjectured, and not yet proved, in \cite{MPT}, where $r=2$.)
Its kernel is also invariant since the target of the morphism is the global sections of local systems, and the projection $p$ is sufficiently general.
This finishes the proofs of Theorems~1 and 2.
\msn
{\bf Remark~2.2.} In \cite[Theorem~0.1]{DS1} we apply the topological cone theorem. This says that $X\cap B_{\ep}$ is identified with the topological cone of the link $X\cap\dd B_{\ep}$ as a stratified space for $\ep$ very very small.
\msn
{\bf 2.3.~Proof of Corollary~1.} Note first that the nearby and vanishing cycle functors commute with the pull-back to a sufficiently general affine space if the space intersects transversally any strata of a Whitney stratification satisfying Thom's $a_f$-condition.
\sk
Let $p:X\to\C^{r-1}$ be a sufficiently general projection. Set
$$U_{\ep,z}:=B_{\ep}\cap p^{-1}(z).$$
Here $B_{\ep}\subset X$ is an open ball of radius $\ep$ with center 0,
and $z\in B'_{\de}\setminus\Si$ in the notation of the proof of Theorem~(2.1). The isomorphism (5) in Theorem~1 says that
$$\Gamma\bl(U_{\ep,z},\F_{(r)}|_{U_{\ep,z}}\br)\q(0<\de\ll\ep\ll 1)
\leqno(2.3.1)$$
is independent of $\ep$ and $z\in B'_{\de}\setminus\Si$. (Note that $C_{\ep}$ in (5) is $Y\cap U_{\ep,z}$.)
Using the Thom-Mather theory, the argument in the proof of Theorem~(2.1) implies that we have the isomorphisms under the restriction morphisms when $\ep$ is slightly changed, and we use fibrations for the case $z$ is changed.
\sk
For $z\in B'_{\de}\setminus\Si$, set
$$W_z=p^{-1}(z)\,(\cong\C^{n-r+2})\subset X,\q f_{W_z}:=f|_{W_z}.$$
We have the isomorphism
$$\F'_{z,(1)}:=\Hc^1(\varphi_{f_{W_z}}A_{W_z})\cong\F_{(r)}|_{W_z},
\leqno(2.3.2)$$
if $z\in B'_{\de}\setminus\Si$ is sufficiently general.
\sk
Since $f$ is a homogeneous polynomial, we have a $\C^*$-action on the sheaf $\F_{(r)}$. We then see that (2.3.1) is independent of any (large) $\ep>0$, where $z\in B'_{\de}\setminus\Si$ is fixed. Note that $p$ is $\C^*$-equivariant (since it is $\C$-linear), and the singular locus of $Y:={\rm Sing}\,f$ is stable by the $\C^*$-action.
Here $z,\de$ are sent to $\al z,\al\de$ respectively by the $\C^*$-action of $\al\in\C^*$.
\sk
We then get the isomorphism
$$\Gamma\bl(W_z,\F_{(r)}|_{W_z}\br)=\rlap{\raise-9pt\h{$\,\,\,\1\scriptstyle\ep$}}\plim\,\Gamma\bl(U_{\ep,z},\F_{(r)}|_{U_{\ep,z}}\br),
\leqno(2.3.3)$$
since the projective system is constant when $\ep\gg 0$.
Here we use the closure of $U_{\ep,z}$ (to get a projective system of sheaves easily) together with the Mittag-Leffler condition as in \cite{Gr}.
(If the Mittag-Leffler condition is not used, we might have to construct a collar neighborhood at infinity for the singular locus of $f_{W_z}$ using the integration of a controlled vector field.
Here the problem is not very difficult, since it can be reduced to the case of local systems using the formulation in Theorem~2 so that the assertion follows from the finite generation of the fundamental group.)
\sk
As for the compatibility with mixed Hodge structure, the isomorphism (10) is factorized as follows:
$$\aligned H^{n-r}(\Ff)&=H^{n-r}(X,\varphi_fA_X)\\ &=H^{n-r}(W,\varphi_{f_W}A_W)\\ &=\Gamma\bl(W,\Hc^{n-r}(\varphi_{f_W}A_W)\br).\endaligned
\leqno(2.3.4)$$
The second and last isomorphisms are induced respectively by the restriction morphism under the inclusion $W\into X$, and the equality (1) together with the Leray type spectral sequence
$$E_2^{p,q}=H^p(W,\Hc^q\varphi_{f_W}A_W)\Longrightarrow H^{p+q}(W,\varphi_{f_W}A_W).$$
The first one is induced by the morphism $\gamma$ in the distinguished triangle
$$(j_0)_!j_0^{-1}\varphi_fA_X\longrightarrow\varphi_fA_X\buildrel{\gamma}\over\longrightarrow(i_0)_*i_0^*\varphi_fA_X\buildrel{+1}\over\longrightarrow,
\leqno(2.3.5)$$
with $i_0:\{0\}\into X$, $\,j_0:X\setminus\{0\}\into X$ natural inclusions.
We get the first isomorphism, since
$$H^j(X,(j_0)_!j_0^{-1}\varphi_fA_X)=0\q(j\in\Z),
\leqno(2.3.6)$$
using the direct image under the projection $\Xt\to\PP^n$ with $\Xt$ the blow-up of $X=\C^{n+1}$ at $0$. This is a line bundle over $\PP^n$, and $X\setminus\{0\}\subset\Xt$ is the associated $\C^*$-bundle. Note that the stratification of $\varphi_fA_X$ is compatible with the $\C^*$-action. The Leray type spectral sequence is compatible with mixed Hodge structure using the ``classical" $t$-structure in \cite[4.6,2]{mhm}. So the isomorphisms in (2.3.4) are compatible with mixed Hodge structure (using \cite[Lemma 4.1]{DMST} for the second isomorphism of (2.3.4)). This finishes the proof of Corollary~1.
\msn
{\bf Remark~2.3.} The second isomorphism of (2.3.4) can be proved by using the {\it weak Lefschetz type theorem\1} which follows from Artin's vanishing theorem \cite{BBD} (see also \cite[2.1.18]{mhp}). We apply this to the direct image of $\varphi_fA_X$ under the compactification $j_X:X\into\PP^{n+1}$. Note that the closure $\Wo$ of $W$ in $\PP^{n+1}$ is the intersection of sufficiently general hyperplanes $H_k\subset\PP^{n+1}$ ($k\in[1,r{-}1]$), and we have the isomorphism
$$\bl(\RR(j_X)_*\varphi_fA_X\br)|_{\Wo}=\RR(j_W)_*\varphi_{f_W}A_W,
\leqno(2.3.7)$$
with $j_W:W\into\Wo$ the natural inclusion, since $W$ is sufficiently general. More precisely, take a Whitney stratification of $\PP^n$ compatible with $\{f\eq 0\}\subset\PP^n$. This implies a stratification of $\PP^{n+1}\setminus\{0\}$ via the projection $\PP^{n+1}\setminus\{0\}\onto\PP^n$ which is a line bundle, where $0\in\PP^{n+1}$ denotes also the image of $0\in X=\C^{n+1}$ by the inclusion $\C^{n+1}\into\PP^{n+1}$. Take the refinement of this stratification such that the zero section $\PP^n=\PP^{n+1}\setminus\C^{n+1}$ is the union of strata of the original stratification. The hyperplane $H_k$ must be transversal to any strata of the stratification of $\mcap_{j<k}H_j$ induced from this stratification inductively.
\msn 
{\bf 2.4.~Reduced hyperplane arrangement case.} Let $f$ be a defining polynomial of a central reduced hyperplane arrangement $\A$ in $X=\C^{n+1}$.
Let $W\subset X$ be a sufficiently general affine subspace of dimension 3.
Set $\A_W:=\A|_W$, $f_W:=f|_W$. Let $C_i$ ($i\in I_{(2)}$) be the connected components of
$$Y_W^1:=Y_{W,\1{\rm sm}}=Y_W\setminus{\rm Sing}^,Y_W\q\h{with}\q Y_W:={\rm Sing}\,f_W=Y\cap W.$$
Let $P_j$ ($j\in I_{(3)})$ be the singular points of $Y_W$.
Note that $I_{(2)},I_{(3)}$ are respectively identified with the sets of edges of $\A_W$ (or equivalently $\A$) with codimension 2 and 3. We write $j<i$ if $P_j\in\Co_i$.
Let $m_i$, $m_j$ be the number of planes of $\A_W$ containing $C_i$ and $P_j$ respectively (that is, the degree of a local defining polynomial). Set $m_{j,i}:=m_j-m_i$ if $j<i$. This is the number of planes of $A_W$ intersecting $\Co_i$ transversally at $P_j$.
Let $e_i$ be the greatest common divisor of the $m_j$ ($j<i$) and $m_i$. It is known that
$$\aligned{\rm rk}\,\Hc^1(\psi_{f_W,\la}\C_W)|_{C_i}&=\begin{cases}m_i\mi 2\pl\delta_{\la,1}&(\la\in\mu_{m_i}),\\ 0&(\la\notin\mu_{m_i}),\end{cases}\\
\dim\Gamma\bl(C_i,\Hc^1(\psi_{f_W,\la}\C_W)|_{C_i}\br)&=\begin{cases}m_i\mi 2\pl\delta_{\la,1}&(\la\in\mu_{e_i}),\\ 0&(\la\notin\mu_{e_i}).\end{cases}\endaligned
\leqno(2.4.1)$$
Indeed, the first assertion can be reduced to the case of Brieskorn-Pham polynomials of 2 variables (using a $\mu$-constant deformation). For the second one, note that the local system monodromy around $P_j\in\Co_i$ is given by $T^{-m_{j,i}}$. (This can be proved by using a point-center blow-up along $P_j$.) It implies that the monodromy group of the local system on $C_i$ is an abelian subgroup of the cyclic group of order $m_i$ generated by the monodromy $T$.
\sk
For $\la\eq 1$, it is quite well-known (see for instance \cite{BS}) that
$$H^j(\Ff)_1=H^j(U)\q\h{with}\q U:=\PP^n\setminus\{f\eq 0\},$$
in particular, $\dim H^1(\Ff)_1\eq d\mi 1$.
\sk
For $i\in I_{(2)}$, $\la\in\C^*\setminus\{1\}$, define
$$I_{(3)}^i:=\{j\in I_{(3)}\mid j<i\},\q I_{(3),\la}:=\{j\in I_{(3)}\mid\la^{m_j}=1\},$$
and similarly for $I_{(2),\la}$ with $m_j$ replaced by $m_i$.
\sk
Using \cite[(3.5.1--2)]{DS1}, the last assertion of (2.4.1) implies that if $H^1(F_{f,0})_{\la}\,{\ne}\,0$, then there is $i\in I_{(2),\la}$ such that $\la^{e_i}=1$, or equivalently, $I_{(3)}^i\subset I_{(3),\la}$. In other words, we get the following.
\msn
{\bf Proposition~2.4.} {\it We have $H^1(F_{f,0})_{\la}\eq 0$, if $I_{(3)}^i\not\subset I_{(3),\la}$ for any $i\in I_{(2),\la}$.}
\ms
We notice, however, that \cite[(3.5.1--2)]{DS1} implies a much stronger assertion as follows.
\msn
{\bf Theorem~2.4.} {\it We have $H^1(F_{f,0})_{\la}\eq 0$, if for any $i\in I_{(2),\la}$, there is $j\in I_{(3)}^i$ such that $H^1(F_{f_W,P_j})_{\la}\eq 0$.}
\ms
Here Corollary~1 is unnecessary and \cite[(3.5.1--2)]{DS1} is sufficient for the vanishing of $H^1(F_{f,0})_{\la}$, since $W$ is identified with a transversal slice to each edge of codimension 3.
In Theorem~(2.4) the hypothesis is much weakened than Proposition~(2.4), since it is rather rare that $H^1(F_{\!g,0})_{\la}\,{\ne}\,0$ for $\la\,{\ne}\,1$ with $g$ a defining polynomial of a central reduced plane arrangement in $\C^3$. Indeed, the first example of such an arrangement has degree 6, and the next ones have 9. Here $\la\eq\exp(\pm2\pi i/3)$, and its dimension can be 1 or 2 when the degree is 9, see for instance \cite[3.1--2]{BDS} and references there.
It seems that $H^1(\Ff)_{\la}\eq 0$ ($\la\,{\ne}\,1$) for central reduced hyperplane arrangements except for some rather special cases, see \cite{Di3}, \cite{MP}, \cite{MPP} for such examples.
Note also that for any homogeneous polynomial $g$, we have $H^{\ssb}(F_{\!g,0})_{\la}\eq 0$ if $\la^{\deg g}\,{\ne}\,1$.
\sk
The above observation shows that we have quite big obstructions almost always if $\la\,{\ne}\,1$, especially when $\la^3\,{\ne}\,1$.
It is then important to get sufficient conditions for the vanishing of $H^1(F_{\!g,0})_{\la}$ ($\la\,{\ne}\,1$) in the 3 variable case. (Here it does not seem easy to use effectively \cite[Corollary 5.3]{Di2}.)
For another type of sufficient condition for the vanishing of $H^1(\Ff)_{\la}$ ($\la\,{\ne}\,1$), see for instance \cite{Ba}.
\ms
Proposition~(2,4) is a rather strong assertion, since it implies a simple proof of the following (see also \cite{ac}).
\msn
{\bf Corollary~2.4} (\cite{BDY}). {\it For a reflection hyperplane arrangement of type $G_{31}$, we have $H^1(\Ff)_{\la}\eq 0$ for $\la\eq\exp(\pm 2\pi i/6)$.}
\msn
{\it Proof.} By Proposition~(2,4) it is sufficient to show the following:
$$\h{For any $C_i\,\,\bl(\1 i\in I_{(2)}\br),\,$ there is $\,P_j\,\,\bl(\1 j\in I_{(3)}^i\br)\,$ with $\,m_j\notin 6\1\Z$,}
\leqno(2.4.2)$$
(forgetting the condition about $m_i$). We will denote by $\Pt_j$, $\Ct_i$ the vector subspaces of $\C^4$ corresponding to $P_j$, $C_i\subset W$ so that $P_j\eq\Pt_j\cap W$, etc. In the notation of \cite[Appendix]{ac}, we may assume that $\Ct_i$ is contained in $\{x_4\eq 0\}$ using the transitivity of the action of the reflection group. We show (2.4.2) by assuming further that $\Pt_j$ is given by the intersection with $\{x_3\eq 0\}$ (exchanging $x_3$ with $x_2$ if $C_i$ is equal to $\{x_3\eq x_4\eq 0\}$). The assertion is then reduced to the following:
$$m_k\notin 6\1\Z\,\,\,\,\,\h{for any}\,\,\,\,P_k\in C_{i_0}.
\leqno(2.4.3)$$
where $C_{i_0}$ corresponds to $\Ct_{i_0}:=\{x_3\eq x_4\eq 0\}$.
Note that $\Ct_i$ corresponding to $C_i$ in (2.4.2) is given by taking the intersection of $\{x_4\eq 0\}$ with a hyperplane of $\A$ defining $\Pt_k\subset\Ct_{i_0}$.
It is then enough to calculate $m_k$ when $\Pt_k$ is the line in $\C^4$ corresponding to $[1\col 0\col 0\col 0]$ or $[1\col 1\col 0\col 0]\in\PP^3$ (using the action of $(\mu_4)^4$ on $\PP^3$). We then always get that $m_k\eq 15$. So the assertion follows. This finishes the proof of Corollary~(2.4).
\msn
{\bf Remark~2.4}\,(i). In the above notation we have more generally $m_j\eq 15,15,4,4,15$ when $\Pt_j$ is the line corresponding to $[1\col 0\col 0\col 0]$, $[1\col 1\col 0\col 0]$, $[2\col 1\col 1\col 0]$, $[1\col 1\col 1\col 0]$, $[1\col 1\col 1\col 1]\in\PP^3$ respectively.
\msn
{\bf Remark~2.4}\,(ii). In \cite{BDY}, \cite{ac}, one considered a hyperplane cut in $\PP^3$, and hence the corresponding vector subspace cut in $\C^4$, but not an affine subspace cut in $\C^4$.
\msn
{\bf Remark~2.4}\,(iii). Consider a reflection arrangement of type $G(m,m,d_X)$ defined by
$$f=\mprod_{1\les j<k\les d_X}\,(x_j^m\mi x_k^m)\q(m\ges 2),$$
where $d_X:=n{+}1$, see \cite[p.\,247 and 280]{OT}. In the case $m\eq 1$, $f$ is non-essential (that is, a polynomial of fewer variables), since $f$ is equivalent to
$$\bl(\mprod_{j=1}^n\,y_j\br)\,\mprod_{1\les j<k\les n}\,(y_j\mi y_k),$$
setting $y_j=x_j\mi x_{n+1}$ ($j\in[1,n]$).
\sk
Let $\zeta=\exp(\pm 2\pi i/3)$. It is known (see \cite{MP}, \cite{MPP}, \cite{Di3}) that
$$\dim H^1(\Ff)_{\zeta}=\begin{cases}2&(d_X\eq 3,\,m\in 3\1\Z),\\1&(d_X\eq 3,\,m\notin 3\1\Z\,\,\,\,\,\h{or}\,\,\,\,\,d_X\eq 4),\\0&(d_X\,{>}\,4).\end{cases}
\leqno(2.4.4)$$
$$H^1(\Ff)_{\la}=0\q\q(\la\notin\mu_3).
\leqno(2.4.5)$$
Using Theorem~(2.4), it is easy to show (2.4.4) for $d_X\,{>}\,4$ and (2.4.5). If $d_X\eq 3$ with $m\notin 3\1\Z$, we can prove (2.4.4--5) using \cite{ESV}, see also \cite{BDS}.
By a calculation using software in \cite{wh}, we can verify (2.4.4--5) for $d_X\eq 3$, $m\,{\les}\,6$.
When $d_X\,{\ges}\,4$, we can prove partially (2.4.4--5) using Remark~(2.4)\,(iv) below, see (3.5) below. Here we see that for any $i\in I_{(2),\zeta}$, we have $H^1(F_{\!f_W,P_j})_{\zeta}\,{\ne}\,0$ for any $j\in I^i_{(3)}$, using (2.4.4) in the case $d_X\eq 4$, $m\notin 3\1\Z$. However, this is not sufficient to show (2.4.4) in this case.
\sk
The reflection arrangement of type $G(m,1,d_X)$ is defined by
$$f=\mprod_{j=1}^{d_X}\,x_j\,\mprod_{1\les j<k\les d_X}\,(x_j^m\mi x_k^m)\q(m\ges 1),$$
see \cite[p.\,279]{OT}. It is known that the following vanishing holds (see \cite{Di3}).
$$H^1(\Ff)_{\la}=0\q(d_X\,{\ges}\,4,\,\,\la\,{\ne}\,1).
\leqno(2.4.6)$$
It is quite easy to show this using Proposition~(2.4), see Remark~(3.5) below.
\msn
{\bf Remark~2.4}\,(iv). In the notation of the introduction, set
$$\F_{(r),\la}:={\rm Ker}(T_s\mi\la)\subset\F_{(r)}\q\h{(with $\C$-coefficients),}$$
so that $\F_{(r),\la,x}=H^1(F_{\!f,x})_{\la}$ for $x\in f^{-1}(0)$. Here we assume that $f$ is reduced and $r\eq n\mi 1$. Since $W\subset\C^{n+1}$ is a sufficiently general affine subspace, we have the canonical isomorphism
$$(\F_{(r),\la})|_W=\Hc^1\psi_{f_W,\la}\C_W.$$
Set
$$I'_{(2),\la}:=\bl\{i\in I_{(2),\la}\,\big|\,\F_{(r),\la,P_j}\,{\ne}\,0\,\,\,\h{for any}\,\,\,j\in I^i_{(3)}\br\}.$$
Assume for any $i\in I'_{(2),\la}$, we have the equality
$${\rm rk}\,(\F_{(r),\la})|_{C_i}\eq 1.$$
We say that $i,i'\in I_{(2),\la}$ are {\it strongly connected,} if there is $P_j\in I^i_{(3)}\cap I^{i'}_{(3)}$ such that
$$\dim\F_{(r),\la,P_j}\eq 1,$$
and moreover the restrictions of a generator $\xi$ of $\F_{(r),\la,P_j}$ to $U\cap C_i$ and $U\cap C_{i'}$ are both non-zero, where $U$ is a sufficiently small open neighborhood of $P_j$ on which $\xi$ is defined.
We can then consider the {\it connected components\1} of $I'_{(2),\la}$.
Here we say that $i,i'\in I'_{(2),\la}$ are connected if there are $i_k\in I'_{(2),\la}$ ($k\in[0,p]$) such that $i_{k-1}$ and $i_k$ are strongly connected for any $k\in[1,p]$ and moreover $i_0\eq i$, $i_p=i'$.
\sk
We say that a connected component is {\it good\1} if it does not contain $i$ such that $i$ and some $i'\in I_{(2),\la}\setminus I'_{(2),\la}$ are strongly connected.
Let $\rho$ be the number of good connected components of $I'_{(2),\la}$.
By \cite[(3.5.1--2)]{DS1}, we then get the inequality
$$\dim H^1(\Ff)_{\la}\les\rho.
\leqno(2.4.7)$$
\msn
{\bf Remark~2.4}\,(v). Let $D$ be a positive integer such that $H^1(F_{\!g,0})_{\la}=0$ for any $\la\notin\mu_3\cup\mu_4$, if $g$ is a defining polynomial of any central reduced plane arrangement in $\C^3$ with $\deg g<D$.
It may be expected that $D$ is at least 15. Note that there is the {\it Hessian\1} arrangement, that is, the reflection arrangement of type $G_{25}$, where $\deg g\eq 12$ and $H^1(F_{\!g,0})_{\la}\,{\ne}\,0$ ($\la\in\mu_4$), see \cite[Ex.\ 6.30 and p.\,232]{OT}, \cite[Ex.\ 3.5]{FY}, \cite[Remark 3.3\,(iii)]{BDS}, and \cite[Theorem~1.7]{Di3}.
Theorem~(2.4) then implies the following: Let $\la\notin\mu_3\cup\mu_4$, and assume for any $i\in I_{(2),\la}$, there is $j\in I^i_{(3)}$ such that $m_j<D$. Then $H^1(\Ff)_{\la}=0$.
\bs\bs
\vbox{\centerline{\bf 3. Examples}
\bsn
In this section we give some example.}
\msn
{\bf 3.1.~Examples with completely vanishing obstructions} ($\,Y^r\eq Y^{r-1}$). Assume that $X=\C^{n+1}$ ($n\ges 3$), and
$$f=g^a\mi h^b\q(a,b\ges 2),$$
where $g,h$ are irreducible homogeneous polynomials with $\deg g\,{\ne}\,\deg h$. Let $Z_g,Z_h\subset\PP^n$ be projective hypersurfaces defined by $g,h$ respectively.
Assume there is a Zariski-closed subset $\Si\subset Z_g\cup Z_h$ of dimension at most $n\mi 3$ such that $(Z_g\cup Z_h)\setminus\Si$ is a divisor with simple normal crossings on $\PP^n\setminus\Si$.
Then $r:=\dim Y=n\mi1$ with $Y:={\rm Sing}\,f$, and in the notation of Proposition~1, we may assume on a neighborhood of $0\in\C^{n+1}$
$$Y^r=Y^{r-1}=C(Z_g\cap Z_h)\setminus C(\Si),$$
where $C(\Si)\subset\C^{n+1}$, etc.\ denotes the affine cone.
The transversal singularity of $f$ at a point of $Y^r$ is locally expressed by
$$y_1^a-y_2^b=0,$$
since $g$, $h$ are part of local coordinates there.
In particular, the transversal Milnor number is $(a\mi 1)(b\mi 1)$.
The global monodromy of the vanishing cohomology sheaf is trivial, since $f$ is the pull-back of $y_1^a-y_2^b$ by the morphism $(g,h):\C^{n+1}\to\C^2$.
\sk
We can show that $Z_g\cap Z_h$ is connected using the Artin's vanishing theorem (see \cite{BBD}) applied to the dual of $\Q_{\PP^n\setminus Z_g}[n]$, $\Q_{Z_g\setminus Z_h}[n{-}1]$. Indeed, it implies the isomorphisms
$$H^0(\PP^n)\simto H^0(Z_g)\simto H^0(Z_g\cap Z_h).$$
We thus get
$${\rm rk}\,H^{n-r}(\Ff)=(a\mi 1)(b\mi 1),$$
\ms
As a special case we have $f=x_0^a\mi h^b$ ($a,b\ges 2$) with $h$ a homogeneous polynomial of variables $x_1,\dots,x_n$ having an isolated singularity at $0\in\C^n$. This generalizes \cite[Examples 6.4 and 6.5]{MPT} using $x^2\pl 2px=(x\pl p)^2\mi p^2$ (where $a\eq b\eq 2$) although Example 6.5 seems to be removed in the third version.
\msn
{\bf Remark~3.1.} If $X$ is non-singular, the shifted constant sheaf $A_X[d_X]$ is self-dual, hence the strong perversity is satisfied for $A_X[d_X]$, and $\F_{(r)}$ is torsion-free, where $r=\dim Y$.
We may assume $A$ is a field to calculate the rank of cohomology groups.
\msn
{\bf 3.2.~Examples with non-vanishing obstructions.} Assume $X=\C^{n+1}$ ($n\ges 3$), and
$$f=(g^a\mi h^b)\1 q\q(a,b\ges 2),$$
with $g,h,q$ irreducible homogeneous polynomials such that $\deg g\,{\ne}\,\deg h$. Let $Z_g,Z_h,Z_q$ be projective hypersurfaces defined by $g,h,q$ respectively.
Assume there is a Zariski-closed subset $\Si\subset Z_g\cup Z_h\cup Z_q$ of dimension at most $n\mi 4$ such that $(Z_g\cup Z_h\cup Z_q)\setminus\Si$ is a divisor with simple normal crossings on $\PP^n\setminus\Si$.
Set $e:={\rm GCD}(a,b)$, $a':=a/e$, $b':=b/e$.
Then $r:=\dim Y=n\mi1$ with $Y:={\rm Sing}\,f$, and we may assume on a neighborhood of $0\in\C^{n+1}$
$$Y^r=C(Z_g\cap Z_h)\cup\mcup_{\la\in\mu_e}\,\bl(Z'_{\la}\cap C(Z_q)\br)\setminus C(\Si'),$$
with
$$Z'_{\la}:=\{g^{a'}\eq\la h^{b'}\}\subset\C^{n+1},\q\Si':=(Z_g\cap Z_h\cap Z_q)\cup\Si.$$
Set
$$\aligned V_0&:=C(Z_g\cap Z_h)\setminus C(\Si'),\q V_{\la}:=Z'_{\la}\cap C(Z_q)\setminus C(\Si'),\\ V_1&:=C(Z_g\cap Z_h\cap Z_q)\setminus C(\Si')\,\bl(=Y^{r-1}\setminus Y^r\br).\endaligned$$
Let $L_0$, $L_{\la}$ ($\la\in\mu_e$) be the restriction of $\F_{(r)}$ to $V_0$ and $V_{\la}$ respectively.
Then the $L_{\la}$ have globally trivial monodromies, and we have
$${\rm rk}\,L_0=(a\mi 1)(b\mi 1),\q{\rm rk}\,L_{\la}=1\,\,\,(\la\in\mu_e).$$
since $f$ has singularity of type $A_1$ at a point of $Z'_{\la}\cap C(Z_q)\setminus C(\Si')$.
\sk
It is known that the unipotent monodromy part of the vanishing cohomology of $x^a\mi y^b$ has rank $e\mi 1$ (since the number of irreducible components is $e$).
By Remark~(3.2) below, we see that the monodromy of $L_0$ around $V_1$ is identified with the inverse of the monodromy $T$, and moreover there is no ``obstruction" for $L_0$ (explained in the introduction) using the long exact sequence there.
However, the latter sequence applied to $f$ implies that there may be non-trivial ``obstructions" for $L_{\la}$, since we get that
$$\aligned{\rm rk}\,\bl((j_0)_*L_0\oplus\mopl_{\la}\,(j_{\la})_*L_{\la}\br)|_{V_1}/(\F_{(r)}|_{V_1})&=\bl((e\mi 1)\pl e\br)\mi\bl((e\mi 1)\pl 1\br)\\&=e\mi 1,\endaligned$$
where $j_0:V_0\into Y\setminus C(\Si')$, $j_{\la}:V_{\la}\into Y\setminus C(\Si')$ are natural inclusions.
\sk
Since $f$ is the pull-back of $(y_1^a\mi y_2^b)\1y_3$ by $(g,h,q):\C^{n+1}\to\C^3$, we have the triviality of global monodromies (including that for the obstruction local system). We then get by \cite[(3.5.2)]{DS1} the equality
$${\rm rk}\,H^{n-r}(\Ff)=\bl((e\mi 1)\pl e\br)\mi(e\mi 1)=e.$$
(We cannot apply Remark~(3.2) below to $f$ at 0 unless $g,h\in\C[x_1,\dots,x_n]$ with $q=x_0$.)
For these rather simple examples, there is no big advantage in cutting $Y$ by general hyperplanes so that $Y$ is replaced by a curve.
\msn
{\bf Remark~3.2.} Assume $f=gx_n\in\C\{x_1,\dots,x_n\}$ with $g\in\C\{x_1,\dots,x_{n-1}\}$. The monodromy of the vanishing cohomology sheaves $\Hc^j\varphi_f\Q_X$ around $x_n=0$ can be identified with the inverse of the monodromy $T$ associated with the vanishing cohomology sheaves $\Hc^j\varphi_g\Q_{X'}$, where $X':=\C^{n-1}$ is identified with $X'{\times}\{c\}$ for $|c|\,{\ne}\,0$ sufficiently small.
(Note that $T$ is the inverse of the Milnor monodromy, see \cite{DS2}.)
Moreover there is a long exact sequence
$$\to H^j(\Ff)\to H^j(F_{\!g,0})\buildrel{\gamma}\over\to H^j(F_{\!g,0})\to H^{j+1}(\Ff)\to,$$
where $\gamma$ is given by $T\mi{\rm id}$.
(This can be shown by using $(\Delta_{\de'}{\times}B'_{\ep'})\cap f^{-1}(\Delta_{\de})\subset\C^n$ with $B'_{\ep'}\subset\C^{n-1}$.)
\sk
This implies that there is no obstruction for the direct image by the inclusion of $\{x_n\,{\ne}\,0\}$, since the latter can be identified with the kernel of $\gamma$.
This may be used for instance to calculate the nearby cycles in the case of a reduced divisor with normal crossings by induction. (In this case, however, we may have to use the action of the symmetric group to determine the morphism to the obstruction.)
\msn
{\bf 3.3.~Another example with non-vanishing obstructions.}
Assume $X=\C^4$, and
$$f=gw\q\h{with}\q g=(x^a\pl y^a)\1 z^{d-a}\pl x^d\pl y^d\in\C[x,y,z],$$
where $d>a\ges 2$. For simplicity, we assume that the projective curve $Z_g:=\{g=0\}\subset\PP^2$ is {\it irreducible}. This is equivalent to that $\{x^a\eq{-}1\}\cap\{x^d\eq{-}1\}\eq\emptyset$, and imposes some conditions on $a,d$ (for instance, $a,d$ cannot be both odd).
It implies that $Z_g$ has only one singular point at $\{x\eq y\eq 0\}$ (using the expression $z^{d-a}=-(x^d+1)/(x^a+1)$ after substituting $y\eq 1$, and similarly with $x,y$ exchanged).
\sk
Let $\Yc$ be the smooth part of $Y:={\rm Sing}\,f$. It has 2 connected components $\Yc_0,\Yc_1$ with
$$\aligned\Yc_0&=\{x\eq y\eq 0\}\setminus(\overline{\Si}_0\cup\overline{\Si}_1),\q\Yc_1=\{g\eq w\eq 0\}\setminus\overline{\Si}_1,\\ \Si_0&=\{x\eq y\eq z\eq 0\}\setminus\{0\},\q\Si_1=\{x\eq y\eq w\eq 0\}\setminus\{0\}.\endaligned$$
Let $L_i$ be the restriction of $\F_{(2)}$ to $\Yc_i$ ($i\eq 0,1$). Then $L_0,L_1$ are local systems of rank $(a\mi 1)^2$ and 1 respectively.
(Note that the transversal singularities along $\Yc_1$ have type $A_1$ in two variables.)
\sk
Since $g$ is a homogeneous polynomial of degree $d$, we see that the local monodromy of $L_0$ around $\Si_0$ is identified with $T^{-d}$, using a well-known relation between the Milnor monodromy and the vertical monodromy in the homogeneous polynomial case.
(This can be shown also by using the blow-up of $\C^4$ along the origin so that the assertion is reduced to a generalization of Remark~(3.2) where $gx_n$ is replaced by $gx_n^d$.)
The local monodromy of $L_0$ around $\Si_1$ is trivial by Remark~(3.2) (where the center of coordinates of $\C^3$ is shifted).
It is easy to see that $L_1$ is constant, since $g$ and $w$ are globally defined (and $\{g=0\}$ intersects $\{w=0\}$ transversally).
Set $e':={\rm GCD}(a,d)$.
We see that
$${\rm rk}\,\Gamma(\Yc_0,L_0)=e'(a\mi 2)\pl 1,\q{\rm rk}\,\Gamma(\Yc_1,L_1)=1.
$$
Indeed, for $h=x^a\pl y^a$, we have
$$\dim H^1(F_{\!h,0},\C)_{\la}=a\mi 2\pl\delta_{\la,1}\q(\la\in\mu_a).$$
\sk
On the other hand, we can show that
$$H^1(F_{\!g,0},\C)_1=0,$$
using the embedded resolution of $Z_g\subset\PP^2$ (since we assume that $Z_g$ is irreducible), see \cite[1.3]{BS}.
Combined with Remark~(3.2), this vanishing implies that there are obstructions for global sections of $L_0$, and these cannot contribute to $H^1(\Ff)$ via \cite[(3.5.2)]{DS1}.
So only the global sections of $L_1$ can contribute to it via {\it loc.\,cit.}
\sk
We apply Remark~(3.2) to a point of $\Si_1$, where $f$ can be expressed as $f':=(x^a\pl y^a)w$ changing coordinates appropriately. We then get that
$${\rm rk}\,H^1(F_{\!f'})=a.$$
Here the contribution from the kernel of $\1T\mi{\rm id}\1$ on $H^1(F_{\!g'})$ has rank $a\mi 1$ (with $g':=x^a\pl y^a)$.
This coincides with the rank of the open direct image of $L_0$ on $\Si_1$.
We then see that there is no obstruction for {\it global\1} sections of $L_1$ using the irreducibility of $Z_g$ (since ${\rm rk}\,L_1=1$), and moreover,
combined with \cite[(3.5.2)]{DS1}, this induces that
$${\rm rk}\,H^1(\Ff)=1.$$
Actually we have to employ the above calculation of $H^1(F_{\!g,0},\C)_1$ together with Remark~(3.2) in order to justify the above assertion.
These calculations are compatible with the ones using software in \cite{wh}, \cite{nwh}, setting $a\eq 2$, $d\eq 4$.
\msn
{\bf 3.4.~Example with partially vanishing obstructions.}
Assume $X=\C^4$, and
$$f\eq x_1^{a_1}x_2^{a_2}\pl y_1^{b_1}y_2^{b_2}\q\q\bl(a_i,b_j\in\Z_{>0}\,\,(i,j\in\{1,2\})\br),$$
where $x_1,x_2,y_1,y_2$ are the coordinates of $\C^4$. Then
$$Y=\{x_1x_2=0\}\cap\{y_1y_2=0\}=\mcup_{i,j\in\{1,2\}}\,\C^2_{x_i,y_j},$$
where $\C^2_{x_i,y_j}\subset\C^4$ is the vector subspace spanned by the $x_i$ and $y_j$-axes ($i,j\in\{1,2\}$).
Set
$$\aligned Y_{i,j}&:=(\C^*)^2_{x_i,y_j}=\{x_iy_j\,{\ne}\,0\}\subset\C^2_{x_i,y_j}\q(i,j\in\{1,2\}),\\
D_{i,j}&:=\begin{cases}\C^*_{x_i}&(j\eq 0,\,\,i\eq 1,2),\\ \C^*_{y_j}&(i\eq 0,\,\,j\eq 1,2).\end{cases}\endaligned$$
We have
$$\aligned Y^2&=\h{$\bigsqcup$}_{i,j\in\{1,2\}}\,Y_{i,j},\\ Y^1\setminus Y^2&=D_{1,0}\sqcup D_{2,0}\sqcup D_{0,1}\sqcup D_{0,2}\\
&=\C^*_{x_1}\sqcup\C^*_{x_2}\sqcup\C^*_{y_1}\sqcup\C^*_{y_2},\endaligned$$
and
$$\Yo_{i,1}\cap\Yo_{i,2}=\Do_{i,0},\q\Yo_{1,j}\cap\Yo_{2,j}=\Do_{0,j}\q(i,j\in\{1,2\}).$$
Cutting these by a general hyperplane $W$, we get a picture as below:
\par
$$\q\h{$\setlength{\unitlength}{5mm}
\begin{picture}(6,6)
\put(1,0){\line(0,1){6}}
\put(0,1){\line(1,0){6}}
\put(5,0){\line(0,1){6}}
\put(0,5){\line(1,0){6}}
\put(-0.5,0.1){$D'_{0,1}$}
\put(5.4,0.1){$D'_{2,0}$}
\put(-0.5,5.4){$D'_{1,0}$}
\put(5.4,5.4){$D'_{0,2}$}
\put(-0.5,2.9){$Y'_{1,1}$}
\put(5.4,2.9){$Y'_{2,2}$}
\put(2.4,5.4){$Y'_{1,2}$}
\put(2.4,.1){$Y'_{2,1}$}
\end{picture}$}$$
\skn
Here $Y'_{i,j}=Y_{i,j}\cap W$, etc.
Note that the 4 lines are in $\C^3$, and are not contained in one plane.
\sk
Let $L_{i,j}$ be the restriction of $\F_{(2)}$ to $Y_{i,j}$ for $i,j\in\{1,2\}$, or to $D_{i,j}$ for $i\eq 0$, $j\eq1,2$ or $j\eq 0$, $i\eq 1,2$.
We denote by $T_{x_i}$, $T_{x_j}$ the monodromy of the local system $L_{i,j}$ around the divisors $\{x_i\eq 0\}$, $\{y_j\eq 0\}$ respectively, if $i,j\in\{1,2\}$. Set
$$a_0:={\rm GCD}(a_1,a_2),\q b_0:={\rm GCD}(b_1,b_2).$$
For $i,j\in\{1,2\}$, let $i',j'\in\{1,2\}$ with $\{i,i'\}=\{1,2\}$, $\{j,j'\}=\{1,2\}$. We have
$$\aligned{\rm rk}\,L_{i,j}=(a_{i'}\mi 1)(b_{j'}\mi 1),&\q{\rm rk}\,L_{i,j}^{T_{x_i}}=(a_0\mi 1)(b_{j'}\mi 1),\\ {\rm rk}\,L_{i,j}^{T_{y_j}}=(a_{i'}\mi 1)(b_0\mi 1),&\q{\rm rk}\,L_{i,j}^{\langle T_{x_i},T_{y_j}\rangle}=(a_0\mi 1)(b_0\mi 1),\endaligned$$
where $\langle T_{x_i},T_{y_j}\rangle$ denotes the subgroup generated by $T_{x_i},T_{y_j}$.
Moreover the first equality holds also for $i\eq i'\eq 0$, $j\in\{1,2\}$ or $j\eq j'\eq 0$, $i\in\{1,2\}$.
These can be verified by applying the Thom-Sebastiani theorem to $x_i^{a_i}x_{i'}^{a_{i'}}\pl y_{j'}^{b_{j'}}$, etc.
Indeed, the vanishing cohomology of $x_{i'}^{a_{i'}}\pl y_{j'}^{b_{j'}}$ with $\C$-coefficients is the tensor product of
$$\mopl_{k=1}^{a_{i'}-1}\,\C u_k,\q\mopl_{l=1}^{b_{j'}-1}\,\C v_l,$$
and the action of the monodromy $T$ is given by $T'\otimes T''$ with
$$T'u_k=\zeta_{a_{i'}}^{\1k}\1u_k,\q T''v_l=\zeta_{b_{j'}}^{\1l}\1v_l,$$
where $\zeta_m:=\exp(2\pi i/m)$ for $m\in\Z_{>0}$.
We see that the local system monodromy $T_{x_i}$ coincides with $(T')^{-a_i}$ (and similarly for $T_{y_j}$).
This can be reduced to the case $f\eq x_1^{a_1}x_2^{a_2}$ using the Thom-Sebastiani theorem. We then get the above assertion on the ranks.
\sk
These calculations show that there is no obstruction as long as one irreducible component $Y_{i,j}$ is considered for $D_{0,j}$ or $D_{i,0}$, although it exists if the direct sum for the two components $Y_{i,j}$, $Y_{i',j}$ or $Y_{i,j}$, $Y_{i,j'}$ is taken, where $i\pl i'\eq j\pl j'\eq 3$. (This can be observed even in the case $f=x_1^{a_1}x_2^{a_2}$ or $y_1^{b_1}y_2^{b_2}$.)
We then get that the {\it obstruction\1} local system on $D_{0,j}$ or $D_{i,0}$ is isomorphic to $L_{0,j}$ or $L_{i,0}$.
By \cite[(3.5.2)]{DS1}, it can be expected that
$${\rm rk}\,H^1(\Ff)=(a_0\mi 1)(b_0\mi 1).$$
However, this is not completely trivial. Strictly speaking, we have to show some {\it compatibility\1} between the restriction morphisms from the global sections of the $L_{i,j}$ for $i,j\in\{1,2\}$ to those of the obstruction local systems which are isomorphic to the $L_{i,j}$ for $i\eq 0$, $j\eq1,2$ or $j\eq 0$, $i\eq 1,2$.
Here each direct factor of the source and target of the morphism (9) in the introduction has the same rank $(a_0\mi 1)(b_0\mi 1)$ (since we take the global sections), but it is not clear whether we have some good generators such that all the restrictions morphisms are simultaneously represented by the identity matrix up to sign (and we would have to see whether the sign is correct).
\sk
Fortunately we can prove the above equality using the Thom-Sebastiani theorem applied to $f\eq x_1^{a_1}x_2^{a_2}\pl y_1^{b_1}y_2^{b_2}$, since ${\rm rk}\,\Ht^0(F_{\!g})\eq a_0\mi 1$ for $g\eq x_1^{a_1}x_2^{a_2}$.
This implies conversely a certain compatibility between the above restriction morphisms.
(In the case $a_i\eq b_j\eq c$ for any $i,j\in\{1,2\}$, it may be possible to verify the compatibility using a group action.)
\msn
{\bf 3.5.~Example of a reflection hyperplane arrangement.} Assume $X=\C^{d_X}$, and a hyperplane arrangement $\A$ is defined by
$$f=\mprod_{1\les j<k\les d_X}\,(x_j^m\mi x_k^m)\q(m\ges 2),$$
where $d_X:=n{+}1$. This is a reflection arrangement of type $G(m,m,d_X)$, see \cite[p.\,247 and 280]{OT}. Setting $\eta:=\exp(2\pi i/m)$, each hyperplane of $\A$ is defined by
$$h_{j,k,p}:=x_j\mi\eta^px_k\q(1\,{\les}\,j\,{<}\,k\,{\les}\,d_X,\,\,p\in\Z/m\Z).$$
\sk
In the notation of (2.4), we have the decomposition
$$I_{(2)}=\msqcup_{a=2}^4\,I^a_{(2)}\q\q\h{with}\q\q I^a_{(2)}:=\bl\{i\in I_{(2)}\,\big|\,|{\rm Supp}\,i|=a\br\}.$$
Here ${\rm Supp}\,i$ is the minimal subset $K_i\subset\{1,\dots,d_X\}$ such that any irreducible component of $\A$ containing $\Ct_i\subset\C^{d_X}$ is defined by a linear combination of the $x_k$ with $\,k\in K_i$. More explicitly, it is equal to $\{j,k\}\cup\{j',k'\}$ if $\Ct_i$ is defined by $h_{j,k,p}\eq h_{j',k',p'}\eq 0$.
We have similarly
$$I_{(3)}=\msqcup_{a=3}^6\,I^a_{(3)}.$$
\sk
For $i\in I^a_{(2)}$, $j\in U^a_{(3)}$, we can verify that
$$\h{$m_i=\begin{cases}\,m&(a\eq 2),\\ \,3&(a\eq 3),\\ \,1{+}1&(a\eq 4),\end{cases}$}\q\q\q\h{$m_j=\begin{cases}\,3m&(a\eq 3),\\ \,6\,\,\,\,\h{or}\,\,\,\,m{+}1&(a\eq 4),\\ \,3{+}1&(a\eq 5),\\ \,1{+}1{+}1&(a\eq 6).\end{cases}$}
\leqno(3.5.1)$$
Here $+$ is used to express the sum of the degrees of {\it indecomposable factors\1} (for instance, $1{+}1$ and $1{+}1{+}1$ mean normal crossing divisors).
There are two cases if $j\in I^4_{(3)}$; the first one is closely related to a remark about the case $m\eq 1$ in Remark~(2.4)\,(iii), and the second one is decomposable. Note that $I^a_{(2)}\eq I^a_{(3)}\eq\emptyset$ if $a\,{>}\,d_X\eq n{+}1$.
\sk
Using Proposition~(2.4) and Remark~(2.4)\,(iv) together with (3.5.1), we can prove partially (2.4.4--5) for $d_X\,{\ges}\,4$ as follows:
$$\aligned\dim H^1(\Ff)_{\la}\les 1&\,\,\,\,\h{if}\,\,\,\,d_X\eq 4,\,\,\la\eq\exp(\pm 2\pi i/3),\\ H^1(\Ff)_{\la}=0&\,\,\,\,\h{if}\,\,\,\,d_X\eq 4,\,\,\la\notin\mu_3\,\,\,\,\,\h{or}\,\,\,\,\,d_X\,{>}\,4.\endaligned
\leqno(3.5.2)$$
(As for the first inequality, the equality actually holds by \cite{MPP}.)
Indeed, the last vanishing follows immediately from (3.5.1) using Proposition~(2.4). The first inequality follows also from (3.5.1) using Remark~(2.4)\,(iv). Here we have to prove that if $C_i\cap C_{i'}\eq\{P_j\}$ for $i,i'\in I^3_{(2)}$, $P_j\in I^4_{(3)}$ with $m_i\eq m_{i'}\eq 3$, $m_j=6$, then $i,i'$ are strongly connected at $j$. This is shown by using an action of a finite group on the central hyperplane arrangement defined by
$$xyz(x\mi y)(y\mi z)(x\mi z)\eq 0.$$
This is equivalent to $(u^2\mi v^2)(v^2\mi w^2)(u^2\mi w^2)\eq 0$, setting $x\eq u\pl v$, $y\eq u\pl w$, $z\eq v\pl w$ (that is, $A_3\eq D_3$), and we can use the natural action of $(\mu_2)^3$ on $\C^3$.
\msn
{\bf Remark~3.5.} A reflection arrangement of type $G(m,1,d_X)$ is defined by
$$f=\mprod_{j=1}^{d_X}\,x_j\,\mprod_{1\les j<k\les d_X}\,(x_j^m\mi x_k^m)\q(m\ges 1),$$
see \cite[p.\,279]{OT}. In this case, we have the following for $i\in I^a_{(2)}$, $j\in U^a_{(3)}$\,:
$$\aligned m_i&=\begin{cases}(m{+}2)&(a\eq 2),\\ \,3\,\,\,\,\h{or}\,\,\,\,1{+}1&(a\eq 3),\\ \,1{+}1&(a\eq 4),\end{cases}\\ m_j&=\begin{cases}\,3(m{+}1)&(a\eq 3),\\ (m{+}2){+}1\,\,\,\,\h{or}\,\,\,\,3{+}1&(a\eq 4),\\ \,3{+}1\,\,\,\,\h{or}\,\,\,\,1{+}1{+}1&(a\eq 5),\\ \,1{+}1{+}1&(a\eq 6).\end{cases}\endaligned
\leqno(3.5.3)$$
Using Proposition~(2.4), this implies that $H^1(\Ff)_{\la}\eq 0$ ($d_X\,{\ges}\,4,\,\la\,{\ne}\,1$), see also \cite{Di1}.
\bs\bs
\vbox{\centerline{\bf Appendix.~Relation between intersection matrices and variations.}
\bsn
The duality for nearby and vanishing cycle functors implies that the determinant of the intersection matrix $S$ coincides with that of $M:={\rm Id}\1{-}\,T$ up to sign in the hypersurface isolated singularity case, where $T$ is the Milnor monodromy, see Remark~(A.1) below.
In the surface $A,D,E$ singularity case, we can verify this by a direct computation as follows.}
\msn
{\bf Case 1} : $A_k$. Let $f(k)$ be the determinant of $-S$ for $A_k$. The intersection matrix $S$ is given by the corresponding Dynkin diagram via the simultaneous resolution as is well known. The matrices $-S$ for $A_4$, $D_5$, $E_6$ are noted below. (These are positive definite.)
Using the inductive expression of determinant, we get the recursive relation
$$f(k)= 2 f(k{-}1) -f(k{-}2),$$
hence $f(k)=k\1{+}\11$. On the other hand, we have
$$\det M=\mprod_{j=1}^k (1\1{-}\1e^{2\pi i\al_j})=(t^{k+1}\1{-}\11)/(t\1{-}\11)|_{t\to 1}=k\1{+}\11,$$
with $\al_i$ the spectral numbers. These are $\dfrac{i}{k{+}1}+1$ for $i\in[1,k]$ in the case of $A_k$, where the equation is given by $x^{k+1}{+}y^2{+}z^2$. So the coincidence follows.
\msn
$$\scalebox{0.8}{$A_4:\begin{pmatrix}2&-1\\-1&2&-1\\&-1&2&-1\\&&-1&2\\\end{pmatrix}$\q$D_5:\begin{pmatrix}2&&-1\\&2&-1\\-1&-1&2&-1\\&&-1&2&-1\\&&&-1&2\end{pmatrix}$\q$E_6:\begin{pmatrix}2&&&-1&\\&2&-1\\&-1&2&-1\\-1&&-1&2&-1\\&&&-1&2&-1\\&&&&-1&2\end{pmatrix}$}$$
\skn
{\bf Case 2} : $D_k$. Let $g(k)$ be the determinant of $-S$ for $D_k$. We have
$$g(k) = 2 f(k{-}1) -2 f(k{-}3) = 2k - 2(k{-}2) = 4,$$
and
$$\det M=(t\1{+}\11)(t^{2k-2}\1{-}\11)/(t^{k-1}\1{-}\11)|_{t\to 1}=4,$$
since the spectral numbers of $D_k$ are $\dfrac{2i{-}1}{2(k{-}1)}+1$ ($i\in[1,k{-}1]$) and $\dfrac{\,3\,}{\,2\,}$, where the equation is given by $x^{k-1}{+}xy^2{+}z^2$.
So we get the coincidence.
\bsn
{\bf Case 3} : $E_k$. Let $h(k)$ be the determinant of $-S$ for $E_k$. We have
$$h(k) = 2 f(k{-}1) - 3 f(k{-}4) = 2k - 3(k{-}3) = 9-k.$$
On the other hand, the spectral numbers $\al_i$ of $E_k$ for $k=6,7,8$ are respectively
$$\aligned\{13,16,19,17,20,23\}/12,\\ \{19,23,25,27,29,31,35\}/18,\\ \{31,37,43,49,41,47,53,59\}/30,\endaligned$$
where the equations are given by $\,x^4{+}y^3{+}z^2$, $\,\,x^3{+}xy^3{+}z^2$, $\,\,x^5{+}y^3{+}z^2\,$ respectively.
\sk
The cyclotomic polynomials $\Phi_{n}(u)$ for $n=12,18,30$ are as follows:
$$\aligned\Psi_{12}(u)&=u^4{-}u^2{+}1,\\ \Psi_{18}(u)&=u^6{-}u^3{+}1,\\ \Psi_{30}(u)&=u^8{+}u^7{-}u^5{-}u^4{-}u^3{+}u{+}1.\endaligned$$
We can calculate the product of $(1\1{-}\1u^{m_k\al_i})$ ($i\in[1,k]$) in $\Q[u]/(\Psi_{m_k}(u))$ using a computer, and get $9\1{-}\1 k$ for $k=6,7,8$ with $m_k=12,18,30$ respectively.
(We cannot replace $\Psi_{m_k}(u)$ for instance with $(u^{m_k}\1{-}\11)/(u\1{-}\11)$.)
This can be done in the case of $E_8$ with Macaulay2, for instance, as follows:
\ms
\vbox{\small\sf\verb#R=QQ[u]; S=R/(u^8+u^7-u^5-u^4-u^3+u^1+1);#
\sk
\verb#(1-u^31)*(1-u^37)*(1-u^43)*(1-u^49)*(1-u^41)*(1-u^47)*(1-u^53)*(1-u^59)#}
\msn
So the coincidence follows.
\msn
{\bf Remark A.1.} Let $f:X\to\Delta$ be a projective morphism from a complex manifold $X$ to a disk $\Delta$ such that $f$ is smooth outside $0\in X$.
There is a well-known relation
$$\langle {\rm can}\,v,w\rangle=\pm\langle v,{\rm var}\,w\rangle,$$
for $v\in H^0(X_0,\psi_f\Z_X[n])$, $w\in H^0(X_0,\varphi_f\Z_X[n])$, concerning the two morphisms
$${\rm can}:\psi_f\Z_X\to\varphi_f\Z_X,\q{\rm var}:\varphi_f\Z_X\to\psi_f\Z_X,$$
such that
$${\rm var}\ssc{\rm can}=T\mi{\rm id},\q{\rm can}\ssc{\rm var}=T\mi{\rm id}.$$
There is no Tate twist, since these morphisms are defined topologically.
(They are {\it different\1} from those used in Hodge theory, see \cite[5.2.3]{mhp}.)
The above relation means that var is identified up to sign with the {\it dual\1} of can, which gives a morphism of homology groups.
\sk
If $v={\rm var}\,v'$, we get the relation
$$\langle {\rm can}\ssc{\rm var}\,v',w\rangle=\pm\langle {\rm var}\,v',{\rm var}\,w\rangle.$$
This implies that if there is a duality for nearby and vanishing cycle functors without up to torsion, then the determinant of the intersection matrix coincides up to sign with that of the variation $T\1{-}\,{\rm id}={\rm can}\ssc{\rm var}$ in the hypersurface isolated singularity case.
(Note that these vanish if the unipotent monodromy part does not vanish.)
\msn
{\bf Remark A.2.} The above relation between the determinants seems to be known to some specialists, see for instance the proof of Proposition 4.7, page 93 in a book of Dimca \cite{Di1}. A formula for the determinant up to sign of the intersection matrices for surface $A,D,E$ singularities is also noted in p.\1\,222 of this book. (One can verify it also by using det in Macaulay2 for instance.)

\end{document}